\documentclass[10pt,twocolumn,english]{IEEEtran}
\usepackage[T1]{fontenc}
\usepackage[latin9]{inputenc}
\usepackage{amsthm}
\usepackage{amsmath}
\usepackage{amssymb}
\usepackage{graphicx}
\usepackage{dsfont}

\usepackage{amsthm,amsmath,amsfonts,amssymb,mathtools}
\usepackage{enumerate}
\usepackage{graphicx}
\usepackage{bbm,cite}
\usepackage{fancyhdr}
\usepackage{marginnote}

\makeatletter
\theoremstyle{plain}

\usepackage{amsthm}
\usepackage{amssymb}

\theoremstyle{plain}
\theoremstyle{remark}
\theoremstyle{plain}
\theoremstyle{plain}

\ifCLASSINFOpdf
\else
\fi
\usepackage{babel}

\usepackage{babel}

\providecommand{\theoremname}{Theorem}

\usepackage{babel}
\providecommand{\theoremname}{Theorem}

\usepackage{babel}
\providecommand{\theoremname}{Theorem}

\usepackage{babel}
\providecommand{\theoremname}{Theorem}
\input{tcilatex}

\makeatother

\usepackage{babel}
\begin{document}

\title{Cloud Radio-Multistatic Radar: Joint Optimization of Code Vector
and Backhaul Quantization}

\author{Shahrouz Khalili, Osvaldo Simeone, \emph{Senior Member, IEEE}, Alexander M. Haimovich,
\emph{Fellow, IEEE}
\thanks{Copyright (c) 2012 IEEE. Personal use of this material is permitted. However, permission to use this material for any other purposes must be obtained from the IEEE by sending a request to pubs-permissions@ieee.org.}\\
\thanks{S. Khalili, O. Simeone and A. M. Haimovich are with CWCSPR, ECE Dept,
NJIT, Newark, USA. E-mail: \{sk669, osvaldo.simeone, haimovic\}@njit.edu.} }
\maketitle
\begin{abstract}
A multistatic radar set-up is considered in which distributed receive
antennas are connected to a Fusion Center (FC) via limited-capacity
backhaul links. Similar to cloud radio access networks in communications,
the receive antennas quantize the received baseband signal before
transmitting it to the FC. The problem of maximizing the detection
performance at the FC jointly over the code vector used by the transmitting
antenna and over the statistics of the noise introduced by backhaul
quantization is investigated. Specifically, adopting the information-theoretic
criterion of the Bhattacharyya distance to evaluate the detection
performance at the FC and information-theoretic measures of the quantization
rate, the problem at hand is addressed via a Block Coordinate Descent
(BCD) method coupled with Majorization-Minimization (MM). Numerical
results demonstrate the advantages of the proposed joint optimization
approach over more conventional solutions that perform separate optimization.
\end{abstract}
\global\long\def\figurename{Fig.}



\section{Introduction}

Waveform design has been a topic of great interest to radar designers,
see, e.g., \cite{Cook2012}, \cite{Rihaczek1969}, \cite{Levanon2004}.
In particular, for the problem of signal detection, the shape of the
transmitted waveform may greatly affect detection performance when
the radar operates in a clutter environment in which detection is
subject to signal-dependent interference. The optimal waveform in
the Neyman-Pearson (NP) sense is studied for monostatic radars in
\cite{DeLong1967}\cite{Kay2007}. With multistatic radars, the NP criterion
affords little insight into optimal waveform design, and information-theoretic
criteria, such as Kullback-Leibler divergence \cite{Kay2009} and
Bhattacharyya distance \cite{batch}, have served in the literature
as tractable alternatives \cite{lam}.

\begin{figure}[t]
\centering \includegraphics[scale=0.3]{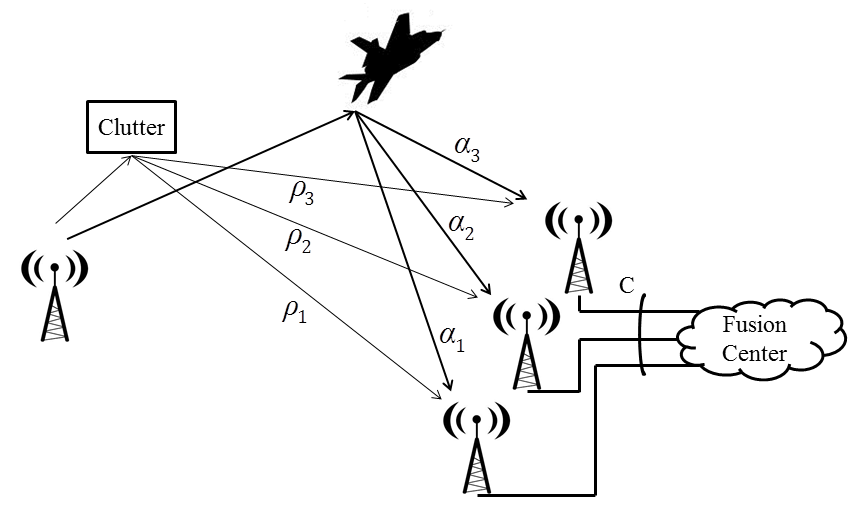}
\protect\caption{Illustration of a CR-MR operating in the presence of a target and
of clutter with $N=3$ receive antennas.}

\label{sys}
\end{figure}

Existing waveform design techniques such as those discussed in \cite{Kay2009,lam},
assume infinite-capacity links between a set of distributed radar
elements and a Fusion Center (FC) that performs target detection (see
Fig. \ref{sys}). In scenarios in which the receive antennas are distributed
over a large geographical area to capture a target's spatial diversity
\cite{him} and no wired backhaul infrastructure is in place, this
assumption should be revised. In fact, in such cases, including deployments
in hostile environments or with moving sensors, the antennas would
be typically connected to the FC through limited-capacity backhaul
links, e.g., microwave radio channels.

In order to cope with the capacity limitations of the backhaul links,
inspired by the cloud radio access architecture in cellular communication
systems \cite{chin}, we assume that the receive sensors quantize
the received baseband signal prior to the transmission to the FC.
Hence, the FC operates on the quantized received baseband signals.
We refer to this system as Cloud Radio-Multistatic Radar (CR-MR).
We formulate and tackle the problem of jointly optimizing over the
code vector and over the operation of the quantizers at the receive
antennas by adopting information-theoretic criteria in Sec. \ref{sec:prob}
and Sec. \ref{sec:solve}, respectively. We observe that, while the
impact of quantization on FC-based sensing systems has been widely
investigated (see, e.g., \cite{ribe}), ours seems to be the first
work to address the joint optimization of code vector and quantization
for multistatic radars. Numerical results are reported in Sec. \ref{sec:num}. 

\textit{Notation}: Bold lowercase letters denote column vectors and
bold uppercase letters denote matrices; $|\mathbf{X}|$ denotes the
determinant of matrix $\mathbf{X}$. $I(X;Y)$ represents the mutual
information between random variables $X$ and $Y$. $\mathcal{CN}(\boldsymbol{\mu},\mathbf{Q})$
is the complex Gaussian distribution with mean vector $\boldsymbol{\mu}$
and covariance matrix $\mathbf{Q}$. $\mathbf{1}$ is a column vector
with all elements equal to one and $[\mathbf{X}]_{m,n}$ denotes the
$(m,n)$ element of matrix $\mathbf{X}$. $\mathds{S}_{+}^{N}$ denotes
the set of symmetric positive semidefinite $N\times N$ matrices.

\section{System Model}

\label{sec:sys} We focus on the CR-MR system shown in Fig. \ref{sys},
in which a transmitter and $N$ receive antennas form a system seeking
to detect the presence of a single stationary target over a clutter
field. The receive antennas are connected to a FC via limited-capacity
backhaul links. While the presented framework is sufficiently general
to accommodate arbitrary backhaul capacity limitations, in this letter,
for simplicity, we adopt the constraint that the total capacity available
for communication between the $N$ receive antennas and the FC is
$C$ bit per received (complex) sample. This scenario captures in
a simple way a backhaul channel that is shared by the receiving antennas.

The radar waveform is a coherent train of $K$ standard pulses with
complex amplitudes forming a code $\mathbf{a}=[a_{1},...,a_{K}]^{T}.$
The pulse repetition intervals are sufficiently large such that the
returns in each pulse interval are due to a single transmitted pulse.
The code design controls the spectral properties of the waveform,
and thus the response of the radar system to the target and clutter.
With respect to each sensor, the target is assumed to obey a Swerling
Type 1 model, i.e., the return has a Rayleigh envelope, which is fixed
over the observation interval. The parameters of the target Rayleigh
envelope are assumed known, and the returns observed by different
sensors are independent. The clutter is assumed homogeneous over the
range of interest, complex-valued Gaussian, with zero-mean and known
variance, fixed over the observations interval and independent between
sensors. Finally, the additive Gaussian noise is assumed to have a
known temporal covariance matrix for each sensor.

Based on the mentioned assumptions, the $K\times1$ discrete-time
signal received by the $n$-th antenna, after matched filtering and
symbol-rate sampling, is given by \cite{lam}
\begin{equation}
\begin{split} & \mathcal{H}_{0}:~\mathbf{r}_{n}=\mathbf{c}_{n}+\mathbf{w}_{n}\\
 & \mathcal{H}_{1}:~\mathbf{r}_{n}=\mathbf{s}_{n}+\mathbf{c}_{n}+\mathbf{w}_{n}~~~n=1,...,N,
\end{split}
\label{rec}
\end{equation}
where the hypotheses $\mathcal{H}_{0}$ and $\mathcal{H}_{1}$ respectively,
represent the absence and presence of a target in a given range resolution
cell; $\mathbf{s}_{n}=\alpha_{n}\mathbf{a}$ is the useful part of
the received signal, with $\alpha_{n}$ being the random complex amplitude
of the target return; $\mathbf{c}_{n}=\rho_{n}\mathbf{a}$ denotes
the clutter, with $\rho_{n}$ being the random clutter complex amplitude;
and $\mathbf{w}_{n}$ is Gaussian noise, accounting for thermal noise,
interference and jamming, which is assumed to be distributed as $\mathcal{CN}(\mathbf{0},\mathbf{M}_{n})$
for some covariance matrix $\mathbf{M}_{n}$. The complex amplitudes
$\alpha_{n}$ and $\rho_{n}$ are independent and distributed as $\mathcal{CN}(0,\sigma_{t,n}^{2})$
and $\mathcal{CN}(0,\sigma_{c,n}^{2})$, respectively. All variables
$\mathbf{w}_{n}$, $\alpha_{n}$ and $\rho_{n}$ are also independent
for different values of $n$. The second order statistics $\sigma_{t,n}^{2}$,
$\sigma_{c,n}^{2}$ and $\mathbf{M}_{n}$ are assumed to be known
to the FC for all $n=1,...,N$, e.g., from prior measurements or prior
information \cite{lam}.

Each receiver quantizes the received vector $\mathbf{r}_{n},$ and
sends the quantized vector $\widetilde{\mathbf{r}}_{n}$ to the FC.
Note that, since the receiver does not know whether the target is
present or not, the quantizer cannot depend on the correct hypothesis
$\mathcal{H}_{0}$ or $\mathcal{H}_{1}$. In order to facilitate analysis
and design, we follow the standard approach of modeling the effect
of quantization by means of an additive quantization noise (see, e.g.,
\cite{gresho}\cite{cover}). The signal received by the FC from the
$n$-th antenna is hence given by
\begin{equation}
\begin{split} & \mathcal{H}_{0}:~\widetilde{\mathbf{r}}_{n}=\mathbf{c}_{n}+\mathbf{w}_{n}+\mathbf{q}_{n}\\
 & \mathcal{H}_{1}:~\widetilde{\mathbf{r}}_{n}=\mathbf{s}_{n}+\mathbf{c}_{n}+\mathbf{w}_{n}+\mathbf{q}_{n},
\end{split}
\label{recq}
\end{equation}
where $\mathbf{q}_{n}\sim\mathcal{CN}(\mathbf{0},\mathbf{Q}_{n})$
is the quantization error vector, which is assumed to be Gaussian
for the sake of tractability. Note that the covariance matrix $\mathbf{Q}_{n}$
defines the shape of the quantization regions and determines the bit
rate required for backhaul communication between antenna $n$ and
the FC \cite{gresho}\cite{cover}.

To set the problem (\ref{recq}) in a standard form, the signal received
at the FC is whitened with respect to the overall additive noise $\mathbf{c}_{n}+\mathbf{w}_{n}+\mathbf{q}_{n}$,
and the returns from all sensors are combined leading to
\begin{equation}
\begin{split} & \mathcal{H}_{0}:~\mathbf{y}\sim\mathcal{CN}(\mathbf{0},\mathbf{I})\\
 & \mathcal{H}_{1}:~\mathbf{y}\sim\mathcal{CN}(\mathbf{0},\mathbf{DSD}+\mathbf{I}),
\end{split}
\label{e:test}
\end{equation}
where $\mathbf{y}=[\mathbf{y}_{1}^{T},...,\mathbf{y}_{N}^{T}]^{T},$
$\mathbf{y}_{n}=\mathbf{D}_{n}\widetilde{\mathbf{r}}_{n},$ $\mathbf{D}_{n}\mathbf{\ }$is
the whitening matrix associated with the $n$-th radar element, $\mathbf{D}_{n}\triangleq(\sigma_{c,n}^{2}\mathbf{aa}^{H}+\mathbf{M}_{n}+\mathbf{Q}_{n})^{-1/2},$
$\mathbf{D}$ is the block diagonal matrix $\mathbf{D}=\limfunc{diag}[\mathbf{D}_{1},...,\mathbf{D}_{N}],$
and $\mathbf{S}$ is the block diagonal matrix $\mathbf{S}=\limfunc{diag}[\sigma_{t,1}^{2}\mathbf{aa}^{H},...,\sigma_{t,N}^{2}\mathbf{aa}^{H}].$
The detection problem described by (\ref{e:test}) has the standard
solution given by the test $\mathbf{y}^{H}\widehat{\mathbf{s}}\gtrless_{\mathcal{H}_{0}}^{\mathcal{H}_{1}}\gamma$,
where $\widehat{\mathbf{s}}=\mathbf{DSD}\left(\mathbf{DSD}+\mathbf{I}\right)^{-1}\mathbf{y}$
and the threshold $\gamma$ is determined from the tolerated false
alarm probability \cite{Kay1998}.

\section{Problem Formulation}

\label{sec:prob} In this section, we aim at finding the optimum code
vector $\mathbf{a}$ and quantization error covariance matrices $\mathbf{Q}_{n}$,
$n=1,...,N$, for a given backhaul capacity constraint $C$. To this
end, for the sake of tractability, we resort to information-theoretic
metrics for both the detection performance and the backhaul capacity
requirements. Specifically, as in \cite{batch} and \cite{lam} (see
also references therein), we adopt the \emph{Bhattacharyya distance}
between the distributions of the quantized received signal (\ref{recq})
under the two hypotheses to evaluate the performance in terms of detection;
moreover, we leverage \emph{rate-distortion theory} to account for
the backhaul capacity requirements \cite{cover}.

For two zero-mean Gaussian distributions with covariance matrix of
$\boldsymbol{\Sigma}_{1}$ and $\boldsymbol{\Sigma}_{2},$ the Bhattacharyya
distance $\mathcal{B}$ is given by $\mathcal{B}=|0.5(\boldsymbol{\Sigma}_{1}+\boldsymbol{\Sigma}_{2})|/\sqrt{|\boldsymbol{\Sigma}_{1}||\boldsymbol{\Sigma}_{2}|}$
\cite{batch}. Therefore, for the signal model (\ref{recq}), the
Bhattacharyya distance can be calculated as $\mathcal{B}=\sum_{n=1}^{N}\mathcal{B}_{n}(\mathbf{a},\mathbf{Q}_{n})$
with \cite{lam}
\begin{equation}
\mathcal{B}_{n}(\mathbf{a},\mathbf{Q}_{n})\triangleq\log\left(\frac{1+0.5\lambda_{n}}{\sqrt{1+\lambda_{n}}}\right),\label{bel_2}
\end{equation}
where we have made explicit the dependence on $\mathbf{a}$ and $\mathbf{Q}_{n},$
and we have defined
\begin{equation}
\lambda_{n}=\sigma_{t,n}^{2}\mathbf{a}^{H}\left(\sigma_{c,n}^{2}\mathbf{a}\mathbf{a}^{H}+\mathbf{M}_{n}+\mathbf{Q}_{n}\right)^{-1}\mathbf{a}.
\end{equation}
We observe that (\ref{bel_2}) is valid under the assumption that
the effect of the quantizers can be well approximated by an additive
Gaussian noise as per (\ref{recq}); in a suitable asymptotic regime, this can be argued by using rate-distortion
theory as briefly discussed below.

The backhaul rate requirement on each $n$th backhaul link is quantified here
by means of the mutual information $I(\mathbf{r}_{n};\widetilde{\mathbf{r}}_{n})$. 
Rate-distortion theory guarantees the existence of
a vector quantizer operating over a large number of  measurement vectors (\ref{rec})
with a rate asymptotically equal to $I(\mathbf{r}_{n};\widetilde{\mathbf{r}}_{n})$
and such that the empirical distribution of the corresponding sequences
($\mathbf{r}_{n};\widetilde{\mathbf{r}}_{n}$) is close to the joint
distribution described by (\ref{recq}) with high probability \cite{cover}.
While the mutual information $I(\mathbf{r}_{n};\widetilde{\mathbf{r}}_{n})$
depends on the actual hypothesis $\mathcal{H}_{0}$ or $\mathcal{H}_{1},$
it is easy to see that $I(\mathbf{r}_{n};\widetilde{\mathbf{r}}_{n})$
is larger under hypothesis $\mathcal{H}_{1}$. Based on this, the mutual
information $I(\mathbf{r}_{n};\widetilde{\mathbf{r}}_{n})$ evaluated
under $\mathcal{H}_{1}$ is adopted here as the measure of the bit
rate required between $n$-th receive antenna and the FC. This can
be easily calculated as $I(\mathbf{r}_{n};\widetilde{\mathbf{r}}_{n})=\mathcal{I}_{n}(\mathbf{a},\mathbf{Q}_{n})$,
with
\begin{equation}
\begin{split}\mathcal{I}_{n}(\textbf{a},\textbf{Q}_{n}) & =\log\left|\textbf{I}+(\textbf{Q}_{n})^{-1}\textbf{M}_{n}\right|\\
 & +\log\left(1+(\sigma_{c,n}^{2}+\sigma_{t,n}^{2})\textbf{a}^{H}(\textbf{Q}_{n}+\textbf{M}_{n})^{-1}\textbf{a}\right),\label{cap_{o}r}
\end{split}
\end{equation}
where again we have made explicit the dependence of mutual information
on $\mathbf{a}$ and $\mathbf{Q}_{n}$.

The problem of maximizing the metric (\ref{bel_2}) over the code
vector $\mathbf{a}$ and the covariance matrices $\mathbf{Q}_{n}$,
for $n=1,...,N$ under total backhaul capacity constraint is stated
as \begin{subequations} \label{opt_pap}
\begin{align}
 & \underset{\mathbf{a},\left\{ \mathbf{Q}_{n}\right\} _{n=1}^{N}}{\mathrm{minimize}}~\sum_{n=1}^{N}\bar{\mathcal{B}}_{n}(\mathbf{a},\mathbf{Q}_{n})\triangleq-\sum_{n=1}^{N}\log\left(\frac{1+0.5\lambda_{n}}{\sqrt{1+\lambda_{n}}}\right)\label{opt_pap_1}\\
 & \mathrm{subject~to~}\lambda_{n}=\sigma_{t,n}^{2}\mathbf{a}^{H}\left(\sigma_{c,n}^{2}\mathbf{a}\mathbf{a}^{H}+\mathbf{M}_{n}+\mathbf{Q}_{n}\right)^{-1}\mathbf{a}\label{opt_pap_2}\\
 & ~~~~~~~~~~~~~\mathbf{Q}_{n}\succeq0~\mathrm{for~all~}n=1,...,N\label{opt_pap_23}\\
 & ~~~~~~~~~~~~~||\mathbf{a}||_{2}^{2}\leq P\label{opt_pap_3}\\
 & ~~~~~~~~~~~~~\sum_{n=1}^{N}\mathcal{I}_{n}(\mathbf{a},\mathbf{Q}_{n})\leq C,\label{opt_pap_4}
\end{align}
where we have formulated the problem as the minimization of the negative
distance $\sum_{n=1}^{N}\bar{\mathcal{B}}_{n}(\mathbf{a},\mathbf{Q}_{n})$,
with $\bar{\mathcal{B}}_{n}(\mathbf{a},\mathbf{Q}_{n})=-\mathcal{B}_{n}(\mathbf{a},\mathbf{Q}_{n})$,
following the standard convention in \cite{boyd}. The power of the
code $\mathbf{a}$ is constrained not to exceed a value $P.$ Note
that the constraint (\ref{opt_pap_4}) ensures that the total transmission
rate between the receive antennas and the FC is smaller than $C$
according to the adopted information-theoretic metrics.

\section{Solution of the Optimization Problem}

\label{sec:solve} The optimization problem in (\ref{opt_pap}) is
not convex, and is hence difficult to solve to obtain a global optimum.
Aiming at obtaining a locally optimal solution, we approach the joint
optimization of $\mathbf{a}$ and $\mathbf{Q}_{n}$, for $n=1,...,N$
in (\ref{opt_pap}) via BCD. Accordingly, at the $m$-th iteration
of the BCD method, the optimum code vector $\mathbf{a}^{(m)}$ is
obtained by solving (\ref{opt_pap}) for matrices $\mathbf{Q}_{n}$
fixed at given values $\mathbf{Q}_{n}^{(m-1)}$ obtained at the previous
iteration; subsequently, the matrices $\mathbf{Q}_{n}^{(m)}$ are
calculated by solving (\ref{opt_pap}) with $\mathbf{a}=\mathbf{a}^{(m)}$.
The steps of BCD algorithm are summarized in Table \ref{al}.
\begin{table}[h]
\protect\caption{Joint Optimization of Code Vector and Quantization Noise Covariances }

\label{al} \centering %
\begin{tabular}{|l|}
\hline
\textbf{Step 0}: Initialize $\textbf{a}^{(0)}\in\mathds{C}^{N}$ and
$\textbf{Q}_{n}^{(0)}\in\mathds{S}_{+}^{N}$ to $n=1,...,N$ for \tabularnewline
~~~~~~~~~feasible values and set $m=1$.\tabularnewline
\textbf{Step 1}: Find $\textbf{a}^{(m)}$ by solving the optimization
problem in (\ref{opt_pap}) when\tabularnewline
~~~~~~~~~$\textbf{Q}_{n}=\textbf{Q}_{n}^{(m-1)}$ via the
MM algorithm (Sec. \ref{sec:a}, eq. (\ref{opt_a_conv})).\tabularnewline
\textbf{Step 2}: Find $\textbf{Q}_{n}^{(m)}$ for $n=1,...,N$ by
solving the optimization \tabularnewline
~~~~~~~~~problem in (\ref{opt_pap}) when $\textbf{a}=\textbf{a}^{(m)}$
via the MM algorithm\tabularnewline
~~~~~~~~~(Sec. \ref{sec:b}, eq. (\ref{opt_q_conv})).\tabularnewline
\textbf{Step 3}: Set $m=m+1$.\tabularnewline
\textbf{Step 4}: Repeat step 1 and step 2 until the convergence is
attained.\tabularnewline
\hline
\end{tabular}
\end{table}

Steps 2 and 3 of Table \ref{al} still require to solve non-convex
problems. Similar to \cite{lam}, we resort to successive convex approximations
by means of the MM technique \cite{mm-conv}. Note that this algorithm,
which combines BCD and MM, coincides with the general-purpose optimization
scheme studied in \cite{iran}. The MM algorithm converges to a local
optimum, and is based on approximating non-convex functions via convex
functions that are locally tight global upper bounds at the current
iterate. Note that in (\ref{opt_pap}) both functions $\bar{\mathcal{B}}_{n}(\mathbf{a},\mathbf{Q}_{n})$
and $\mathcal{I}_{n}(\mathbf{a},\mathbf{Q}_{n})$ are non-convex in
$\mathbf{a}$ and $\mathbf{Q}_{n}$. Given a non-convex function $f(\mathbf{x})$
of a generic variable $\mathbf{x}$, the MM algorithm at the $l$-th
iteration substitutes the function $f(\mathbf{x})$ with a convex
approximation $f(\mathbf{x}|\mathbf{x}^{[l-1]})$ of $f(\mathbf{x})$
at the current solution $\mathbf{x}^{[l-1]}$ that satisfies the global
upper bound property \end{subequations}
\begin{equation}
f(\mathbf{x}|\mathbf{x}^{[l-1]})\geq f(\mathbf{x}),\label{mm_ex}
\end{equation}
for all $\mathbf{x}$ in the domain, along with the local tightness
condition
\begin{equation}
f(\mathbf{x}^{[l-1]}|\mathbf{x}^{[l-1]})=f(\mathbf{x}^{[l-1]}).\label{mm_ex_2}
\end{equation}
properties guarantee the feasibility of all iterates and convergence
to a local optimum \cite{mm-conv}. We emphasize that we are using
superscript ${(m)}$ to identify the iterations of the outer loop
described by Table \ref{al}, and the superscript ${[l]}$ as the
index of the inner iteration of the MM algorithm. In Sec. \ref{sec:a}
and Sec. \ref{sec:b}, we discuss the application of the MM algorithm
to perform Step 1 and Step 2 in Table \ref{al}.

\subsection{Step 1}

\label{sec:a} At Step 1, the goal is to obtain the optimal value
of $\mathbf{a}^{(m)}$ for problem (\ref{opt_pap}) given $\mathbf{Q}_{n}=\mathbf{Q}_{n}^{(m-1)}$
for all $n=1,...,N$. To this end, we apply the MM algorithm as follows.
A convex locally tight upper bound $\bar{\mathcal{B}}(\mathbf{a},\mathbf{Q}_{n}|\mathbf{a}^{[l]})$
of $\bar{\mathcal{B}}_{n}(\mathbf{a},\mathbf{Q}_{n})$ was derived
in \cite{lam} and is given by
\begin{equation}
\bar{\mathcal{B}}_{n}(\textbf{a},\textbf{Q}_{n}|\textbf{a}^{[l]})=\phi_{n}^{[l]}\textbf{a}^{H}\left((\textbf{M}_{n}+\textbf{Q}_{n})^{-1}\right)\textbf{a}-\mathfrak{R}\left(\left(\textbf{d}_{n}^{[l]}\right)^{H}\textbf{a}\right),
\end{equation}
where $\gamma_{n}=\sigma_{t,n}^{2}/\sigma_{c,n}^{2},~\beta_{n}=\sigma_{c,n}^{2},~\lambda_{n}^{[l]}=\gamma_{n}-\gamma_{n}/(1+\beta_{n}(\mathbf{a}^{[l]})^{H}(\mathbf{M}_{n}+\mathbf{Q}_{n})^{-1}\mathbf{a}^{[l]})$
and
\begin{equation}
\begin{split} & \phi_{n}^{[l]}=\frac{\beta_{n}}{1+\beta_{n}y_{n}^{[l]}}+\beta_{n}(1+0.5\gamma_{n})+\frac{0.5\gamma_{n}}{1+\lambda_{n}^{[l]}}\frac{\beta_{n}}{(1+\beta_{n}y_{n}^{[l]})^{2}}\\
 & \textbf{d}_{n}^{[l]}=\frac{2\beta_{n}(1+0.5\gamma_{n})}{1+\beta_{n}y_{n}^{[l]}(1+0.5\gamma_{n})}\\
 & ~~~~+2\beta_{n}(1+0.5\gamma_{n})(\textbf{M}_{n}+\textbf{Q}_{n})^{-1}\textbf{a}^{[l]}\\
 & x^{[l]}=(\sigma_{c,n}^{2}+\sigma_{t,n}^{2})(\textbf{a}^{[l]})^{H}(\textbf{Q}_{n}+\textbf{M}_{n})^{-1}\textbf{a}^{[l]},
\end{split}
\end{equation}
with $\mathbf{a}^{[l]}$ being the value of $\mathbf{a}$ obtained
at the $l$-th iteration of the MM algorithm and $y_{n}^{[l]}=(\mathbf{a}^{[l]})^{H}(\mathbf{M}_{n}+\mathbf{Q}_{n})^{-1}\mathbf{a}^{[l]}$.
A bound with the desired property can also be easily derived for $\mathcal{I}_{n}(\mathbf{a},\mathbf{Q}_{n})$
by using the inequality $\log\left(1+x\right)\leqslant\log(1+x^{[l]})+1/(1+x^{[l]})(x-x^{[l]})$,
for $x=(\sigma_{c,n}^{2}+\sigma_{t,n}^{2})\textbf{a}^{H}(\textbf{Q}_{n}+\textbf{M}_{n})^{-1}\textbf{a}$,
leading to
\begin{equation}
\begin{split}\mathcal{I}_{n}(\textbf{a},\textbf{Q}_{n}| & \textbf{a}^{[l]})=\log\left|\textbf{I}+(\textbf{Q}_{n})^{-1}\textbf{M}_{n}\right|+\log(1+x^{[l]})\\
 & +\frac{1}{1+x^{[l]}}\big((\sigma_{c,n}^{2}+\sigma_{t,n}^{2})\textbf{a}^{H}(\textbf{Q}_{n}+\textbf{M}_{n})^{-1}\textbf{a}-x^{[l]}\big).
\end{split}
\label{cap_1}
\end{equation}
MM algorithm then prescribes the solution of the following convex
optimization problem iteratively, until convergence is attained:

\begin{subequations} \label{opt_a_conv}
\begin{align}
\mathbf{a}^{[l]}= & \mathrm{arg}~\underset{\mathbf{a}}{\min}~\sum_{n=1}^{N}\bar{\mathcal{B}}_{n}(\mathbf{a},\mathbf{Q}_{n}^{(m-1)}|\mathbf{a}^{[l-1]})\\
 & \mathrm{subject~to~}||\mathbf{a}||_{2}^{2}\leq P\\
 & ~~~~~~~~~~~~~\sum_{n=1}^{N}\mathcal{I}_{n}(\mathbf{a},\mathbf{Q}_{n}^{(m-1)}|\mathbf{a}^{[l-1]})\leqslant C.
\end{align}
\end{subequations}

\subsection{Step 2}

\label{sec:b} At Step 2, the matrices $\mathbf{Q}_{n}^{(m)}$ are
obtained for a given $\mathbf{a}=\mathbf{a}^{(m)}$. Similar to (\ref{cap_1}),
upper bounds with the desired properties are derived for functions
$\mathcal{I}_{n}(\mathbf{a},\mathbf{Q}_{n})$ and $\bar{\mathcal{B}}_{n}(\mathbf{a},\mathbf{Q}_{n})$
as follows

\begin{equation}
\begin{split}\mathcal{I}_{n}(\textbf{a}, & \textbf{Q}_{n}|\textbf{Q}_{n}^{[l]})=\log|\textbf{Q}_{n}^{[l]}+(\sigma_{c,n}^{2}+\sigma_{t,n}^{2})\textbf{a}\textbf{a}^{H}+\textbf{M}_{n}|\\
 & -\log|\textbf{Q}_{n}|+\sum_{k=1}^{N}\mathrm{Tr}\Big(\big(\textbf{Q}_{n}^{[l]}+(\sigma_{c,n}^{2}+\sigma_{t,n}^{2})\textbf{a}\textbf{a}^{H}\\
 & +\textbf{M}_{n}\big)^{-1}(\textbf{Q}_{n}-\textbf{Q}_{n}^{[l]})\Big)
\end{split}
\label{drive_q_cons}
\end{equation}
and
\begin{equation}
\begin{split}\bar{\mathcal{B}}_{n}(\textbf{a}, & \textbf{Q}_{n}|\textbf{Q}_{n}^{[l]})=-\log|(\sigma_{c,n}^{2}+0.5\sigma_{t,n}^{2})\textbf{a}\textbf{a}^{H}+\textbf{Q}_{n}+\textbf{M}_{n}|\\
 & +0.5\mathrm{Tr}\left(\Big(\big(\sigma_{c,n}^{2}\textbf{a}\textbf{a}^{H}+\textbf{M}_{n}+\textbf{Q}_{n}^{[l]}\big)^{-1}\Big)^{T}\textbf{Q}_{n}\right)\\
 & +0.5\mathrm{Tr}\left(\Big(\big((\sigma_{c,n}^{2}+\sigma_{t,n}^{2})\textbf{a}\textbf{a}^{H}+\textbf{M}_{n}+\textbf{Q}_{n}^{[l]}\big)^{-1}\Big)^{T}\textbf{Q}_{n}\right).
\end{split}
\label{drive_q_cons_2}
\end{equation}
The MM algorithm then evaluates the matrices $\mathbf{Q}^{(m)}=[\mathbf{Q}_{1}^{(m)},...,\mathbf{Q}_{N}^{(m)}]$
by solving the following convex optimization problem iteratively,
until convergence is attained: \begin{subequations} \label{opt_q_conv}
\begin{align}
 & \mathbf{Q}^{[l]}=\mathrm{arg}~\underset{\mathbf{Q}}{\min}~\sum_{n=1}^{N}\bar{\mathcal{B}}_{n}(\mathbf{a}^{(m)},\mathbf{Q}_{n}|\mathbf{Q}_{n}^{[l-1]})\\
 & \mathrm{subject~to~}\sum_{n=1}^{N}\mathcal{I}_{n}(\mathbf{a}^{(m)},\mathbf{Q}_{n}|\mathbf{Q}_{n}^{[l-1]})\leqslant C\\
 & ~~~~~~~~~~~~~~\mathbf{Q}_{n}\succeq0~\mathrm{for~all~}n=1,...,N.
\end{align}
\end{subequations}

\section{Numerical Results and Conclusion Remarks}

\label{sec:num} In this section, the performance of the proposed
algorithm that performs joint optimization of the code vector $\mathbf{a}$
and of the quantization noise covariance matrices $\mathbf{Q}_{n}$
for $n=1,...,N$ is investigated via numerical results. For reference,
we consider the performance of the following strategies: (\textit{i})
\textit{No optimization}: Set $\mathbf{a}=\sqrt{P/K}\mathbf{1}$ and
$\mathbf{Q}_{n}=\epsilon\mathbf{I}$, for $n=1,...,N$, where $\epsilon$
is a constant that is found by satisfying the constraint (\ref{opt_pap_4})
with equality; (\textit{ii}) \textit{Code vector optimization}: Optimize
the code vector $\mathbf{a}$ by using the algorithm in \cite{lam},
which is given in Table \ref{al} by setting $\mathbf{Q}_{n}=0$ for
$n=1,...,N$, and set $\mathbf{Q}_{n}=\epsilon\mathbf{I}$, for $n=1,...,N$,
as explained above. 
\textit{Quantization noise optimization}: Set $\mathbf{a}=\sqrt{P/K}\mathbf{1}$
and optimize the covariance matrices $\mathbf{Q}_{n}$ as per Step
2 in Table \ref{al}. (\textit{iv}): \textit{Joint optimization of
code vector and quantization noise}: The code vector $\mathbf{a}$
and the covariance matrices $\mathbf{Q}_{n}$ are optimized jointly
by using the algorithm in Table \ref{al}. Throughout, we set the
number of receive antennas as $N=3$, the length of the code vector
as $K=6$, the variance of the target amplitudes as $\sigma_{t,n}^{2}=1$,
for $n=1,...,N$, the variance of the clutter amplitudes as $\sigma_{c,1}^{2}=0.125,~\sigma_{c,2}^{2}=0.25$,
and $\sigma_{c,3}^{2}=0.5$, the transmitted power as $P=10$ dB and
the noise covariance matrices as $[\mathbf{M}_{n}]_{m,k}=(1-0.15n)^{|m-k|}$
as in \cite{lam}.
\begin{figure}[t]
\centering \includegraphics[scale=0.5]{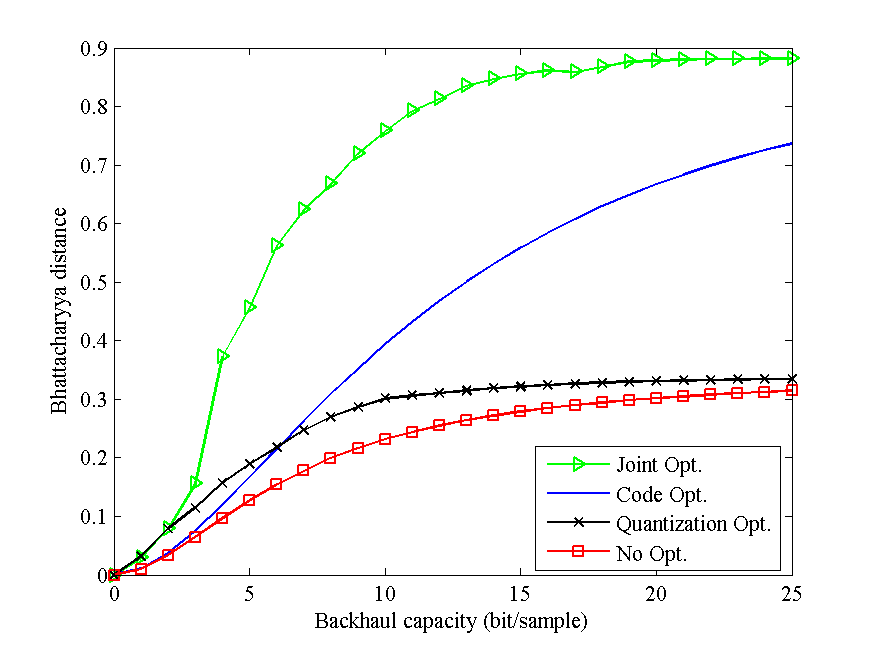}
\protect\caption{Bhattacharyya distance vs. backhaul capacity with $N=3$, $K=6$,
$\sigma_{t,n}^{2}=1$, for $n=1,...,N$, $\sigma_{c,1}^{2}=0.125,~\sigma_{c,2}^{2}=0.25,~\sigma_{c,3}^{2}=0.5$
and $[\mathbf{M}_{n}]_{m,k}=(1-0.15n)^{|m-k|}$.}

\label{bah}
\end{figure}

\begin{figure}[t]
\centering \includegraphics[scale=0.5]{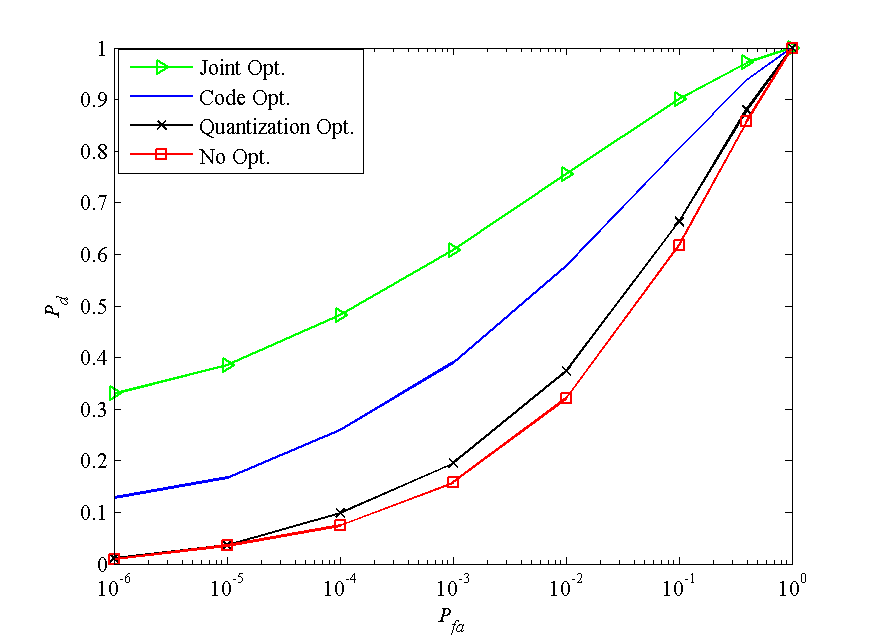} 
\protect\caption{ROC curves with $N=3$, $K=6$, $\sigma_{t,n}^{2}=1$, for $n=1,...,N$,
$\sigma_{c,1}^{2}=0.125,~\sigma_{c,2}^{2}=0.25,~\sigma_{c,3}^{2}=0.5$
and $[\mathbf{M}_{n}]_{m,k}=(1-0.15n)^{|m-k|}$ and $C=15$ bit/sample.}

\label{roc}
\end{figure}

In Fig. \ref{bah} the Bhattacharyya distance is plotted versus the
available backhaul capacity $C$. For sufficiently large values of
$C$, optimizing the code vector only has significant gains as discussed
in \cite{lam}. However, for intermediate values of $C$, it is more
advantageous to properly design the quantization noise. The proposed
joint optimization of code vector and quantization noise is seen to
be beneficial over the separate optimization strategies across all
values of $C$.

Fig. \ref{roc} plots the Receiving Operating Characteristic (ROC),
i.e., the false alarm probability $P_{fa}$ versus the detection probability
$P_{d}$, for $C=15$ bit/sample. The curve was evaluated via Monte
Carlo simulations by implementing the optimum test detector \cite{lam}.
It is seen that the proposed joint optimization method provides remarkable
gains, while, given the backhaul limitations, optimizing only the
code vector leads to significantly smaller advantages. 
 \bibliographystyle{IEEEtran}
\bibliography{refrences}


\end{document}